\documentclass[a4paper, 12pt] {article}
\usepackage[cp1251]{inputenc}
\usepackage[russian, english]{babel}
\usepackage{amsfonts}
\usepackage{mathtext}
\usepackage{cite}

\usepackage{graphicx}

\usepackage{verbatim}
\usepackage{amsmath}
\usepackage{amsthm}
\usepackage{amssymb}
\usepackage{delarray}
\usepackage{mathrsfs}

\begin{document}

\thispagestyle{empty}

\begin{center}
{\bf UNIFORM STABILITY OF THE INVERSE SPECTRAL PROBLEM FOR A CONVOLUTION INTEGRO-DIFFERENTIAL OPERATOR}
\end{center}

\begin{center}
{\large\bf Sergey Buterin}\footnote{Department of Mathematics, Saratov State University, {\it email: buterinsa@info.sgu.ru}}
\end{center}

{\bf Abstract.} The operator of double differentiation, perturbed by the composition of the differentiation operator and a convolution one,
on a finite interval with Dirichlet boundary conditions is considered. We obtain uniform stability of recovering the convolution kernel
from the spectrum in a weighted $L_2$-norm and in a weighted uniform norm. For this purpose, we successively prove uniform stability of
each step of the algorithm for solving this inverse problem in both the norms. Besides justifying the numerical computations, the obtained
results reveal some essential difference from the classical inverse Sturm--Liouville problem.


\\

{\large\bf 1. Introduction}
\\

In this paper, we establish stability of the inverse spectral problem for one important and illustrative class of integro-differential
operators. As is known, this property is vital for justifying numerical algorithms and has a local nature, since it guaranties that
''small'' deviations of any fixed input data caused by measurement or truncation errors may lead only to ''small'' deviations of the
solution. However, the target type of stability is stronger than the usual {\it local} one and belongs to the so-called {\it uniform}
stability, which involves uniform estimates.

Inverse spectral problems consist in recovering operators from their spectral characteristics. The most complete results in the inverse
spectral theory are known for differential operators (see monographs \cite{Mar77,Lev84,FY01,Yur02} and references therein). The first
substantial study in this direction (after the pioneering work of Ambarzumian \cite{A} recently having given a 90-year anniversary to this
topic) was carried out by Borg \cite{Bor}, who proved that the real-valued potential $q(x)\in L_2(0,\pi)$ in the Sturm--Liouville equation
\begin{equation}\label{SLeq}
-y''+q(x)y=\lambda y, \quad 0<x<\pi,
\end{equation}
is uniquely determined by specifying the spectra $\{\lambda_{k,j}\},\;j=0,1,$ of two boundary value problems
${\cal L}_j( q),\;j=0,1,$ for equation (\ref{SLeq}) with one common boundary condition, say:
$$
y(0)=y^{(j)}(\pi)=0,
$$
respectively. For complex-valued potentials, i.e. in the non-selfadjoint case, this uniqueness result was generalized by Karaseva
\cite{Kar}. It is well-known that the following asymptotics holds:
\begin{equation}\label{asympt}
\lambda_{k,j}=\rho_{k,j}^2, \quad \rho_{k,j}=k-\frac{j}{2}+\frac\omega{2k}+\frac{\varkappa_{k,j}}k, \quad \{\varkappa_{k,j}\}\in l_2, \quad
k\ge1, \quad j=0,1.
\end{equation}
Here $\omega=\frac1\pi\int_0^\pi q(x)\,dx.$ Borg \cite{Bor} also established local solvability and local stability of the corresponding
inverse problem. Specifically, the following theorem holds (see also \cite{FY01}).

\medskip
{\bf Theorem 1. }{\it For any model real-valued potential $q(x)\in L_2(0,\pi)$ there exists $\delta>0$ such that if arbitrary real
sequences $\{\tilde \lambda_{k,j}\}_{k\ge1},\;j=0,1,$ satisfy the condition
$$
\Omega:=\sqrt{\sum_{k=1}^\infty\Big(|\lambda_{k,0}-\tilde\lambda_{k,0}|^2+|\lambda_{k,1}-\tilde\lambda_{k,1}|^2\Big)} \le\delta,
$$
then there exists a unique function $\tilde q(x)\in L_2(0,\pi)$ such that $\{\tilde \lambda_{k,j}\}_{k\ge1},\,j=0,1,$ are the spectra of
the problems ${\cal L}_j(\tilde q),\;j=0,1,$ respectively. Moreover,
\begin{equation}\label{est_q}
\|q-\tilde q\|_{L_2(0,\pi)}\le C\Omega,
\end{equation}
where $C$ depends only on $q(x).$}

\medskip
The original proof of Theorem~1 is applicable also for complex-valued potentials, but with the requirement of simplicity of the spectra. In
\cite{ButKuz}, Theorem~1 was generalized for arbitrary multiple spectra, i.e. it remains completely true after replacing all entries of
''real'' with ''complex''. However, in the self-adjoint case or, equivalently, when the function $q(x)$ is real-valued, Marchenko and
Ostrovskii \cite{MO} proved {\it global} solvability of this inverse problem stated in the following theorem.

\medskip
{\bf Theorem 2. }{\it Two sequences $\{\lambda_{k,0}\}_{k\ge1}$ and $\{\lambda_{k,1}\}_{k\ge1}$ are the spectra of the boundary value
problems ${\cal L}_0(q)$ and ${\cal L}_1(q),$ respectively, with a common real-valued potential $q(x)\in L_2(0,\pi)$ if and only if they
are real, have the asymptotics (\ref{asympt}) with a common $\omega$ and interlace in the following way:
$\lambda_{k,1}<\lambda_{k,0}<\lambda_{k+1,1},\;k\ge1.$}

\medskip
Global solvability of the inverse problem inspires one to request its {\it uniform} stability, i.e. when in estimate (\ref{est_q}) the
potential $q(x)$ is also not fixed. The question of uniform stability was raised by Savchuk and Shkalikov in \cite{SavShk}, where they
established, in particular, that for unfixed both $q(x)$ and $\tilde q(x)$ the constant $C$ in (\ref{est_q}) can increase only due to the
following causes:

I. Deviations of the spectra from those corresponding to the zero potential become large;

II. Square roots of neighboring eigenvalues in any pair of problems become too close.\\
Thus, the following refinement of the stability part in Theorem~1 follows from results of \cite{SavShk}.

\medskip
{\bf Theorem 3. }{\it For any $r>0$ and $h\in(0,1/2)$ there exists $C=C(r,h)>0$ (i.e. depending only on $r$ and $h\!$) such that estimate
(\ref{est_q}) holds as soon as:

(I) $\|\{\varkappa_{k,j}\}\|_{l_2}\le r$ and $\|\{\tilde\varkappa_{k,j}\}\|_{l_2}\le r$ for $j=0,1$ and $|\omega|\le r;$ as well as

(II) $h\le \rho_{k,1}+h\le\rho_{k,0}\le\rho_{k+1,1}-h$ and $h\le \tilde\rho_{k,1}+h\le\tilde\rho_{k,0}\le\tilde\rho_{k+1,1}-h$ for
$k\ge1.$\\
Here $\tilde\rho_{k,j}$ and $\tilde\varkappa_{k,j}$ are analogous to $\rho_{k,j}$ and $\varkappa_{k,j},$ respectively, but for the
potential $\tilde q(x).$}

\medskip
Note that $\tilde\omega\ne\omega$ implies $\Omega=\infty.$ In \cite{SavShk}, a different metric is used, which, in particular, admits
different mean values $\omega$ and $\tilde\omega.$ Moreover, in \cite{SavShk}, a uniform two-sided generalization of estimate (\ref{est_q})
was obtained for potentials from a continuous scale of Sobolev spaces with different smoothness indices. Further aspects of uniform
stability for the inverse Sturm--Liouville problem were studied in \cite{SavShk-14}, while local stability in the uniform norm was obtained
in \cite{FY01}.

Restriction (II) in Theorem~3 can be partially explained by the fact that unbounded rapprochement of eigenvalues related to different
problems closes the situation when the inverse problem loses its solvability. We note that for uniform stability of the inverse problem for
the integro-differential operator considered below no similar restriction is necessary.

For integro-differential operators and other classes of nonlocal ones, the classical methods that give global solution of inverse problems
for differential operators (the transformation operator method \cite{Mar77,Lev84,FY01} and the method of spectral mappings
\cite{FY01,Yur02}) do not work. Various aspects of the inverse spectral theory for integro-differential operators were studied in
\cite{Mal, Yur84, Er, Yur91, But04, But07, Kur, But10, KurSh10, W,Yur14, BCh15,BS,BB17,Yur17-1, BondBut-18,
But18,But18-2,Bon18-1,ButVas18,Ign18,Ign18-2,Bon18-2, Zol18, Bon19-1, Bon19-2, Bon19-3,But19} and other works. In particular, in
\cite{Yur91}, developing the idea of Borg's method, Yurko proved local solvability and local stability for the inverse problem of
recovering a convolutional perturbation of the Sturm--Liouville operator from the spectrum. Later on, in \cite{But10} the author proved
global solvability of this inverse problem by another method developing the approach suggested in \cite{But07}, where, in turn, a boundary
value problem ${\mathcal L}={\mathcal L}(M)$ of the form
\begin{equation}\label{eq1}
\ell y:=-y''+ \int_0^x M(x-t)y'(t)\,dt=\lambda y, \;\; 0<x<\pi, \quad y(0)=y(\pi)=0
\end{equation}
was considered with a complex-valued function $M(x)\in L_{2,\pi}:=\{f(x):(\pi-x)f(x)\in L_2(0,\pi)\},$ and the following inverse problem
was studied:

\medskip
{\bf Inverse Problem 1.} Given the spectrum $\{\lambda_k\}$ of ${\mathcal L},$ find the function $M(x).$

\medskip
Here we obtain uniform stability of Inverse Problem~1. Note that even local stability of this inverse problem does not follow from results
of \cite{Yur91}. The operator $\ell$ is quite illustrative, because any results obtained for $\ell$ can usually be generalized to more
complicated classes of integro-differential operators. The following theorem is a fusion of Theorems~1.1 and~1.3 in \cite{But07}.

\medskip
{\bf Theorem 4. }{\it (i) For an arbitrary sequence of complex numbers $\{\lambda_k\}$ to be the spectrum of a boundary value problem
${\mathcal L}$ (which, in turn, is determined by its spectrum uniquely) it is necessary and sufficient to have the asymptotics
\begin{equation}\label{eq3}
\lambda_k=\rho_k^2, \quad \rho_k=k+\varkappa_k, \quad \{\varkappa_k\}\in l_2, \quad k\ge1.
\end{equation}

(ii) The function $M(x)$ satisfies the additional smoothness condition: $M(x)\in W^1_2[0,T]$ for each $T\in(0,\pi)$ and $M'(x)\in
L_{2,\pi}$ if and only if
$$
\lambda_k=\Big(k+\frac{A}{k}+\frac{\varkappa_{k,1}}k\Big)^2, \quad \{\varkappa_{k,1}\}\in l_2, \quad A - const.
$$
Moreover, $M(0)=2A.$}

\medskip
Furthermore, the function $M(x)$ can be found by the following algorithm.

\medskip
{\bf Algorithm 1.} {\it Let the spectrum $\{\lambda_k\}$ of a certain problem ${\mathcal L}(M)$ be given.

1) Construct the function $w(x)$ using the formulae
\begin{equation}\label{ser w}
w(x)=\frac2\pi\sum_{k=1}^\infty k\Delta(k^2)\sin kx, \quad \Delta(\lambda)=\pi\prod_{k=1}^\infty\frac{\lambda_k-\lambda}{k^2};
\end{equation}

2) Find the function $N(x)$ by solving the nonlinear integral equation
\begin{equation}\label{main_eq}
w(\pi-x)=\sum_{\nu=1}^\infty \frac{(\pi-x)^\nu}{\nu!} N^{*\nu}(x),
\end{equation}
where $f^{*1}(x)=f(x)$ and $f^{*(\nu+1)}(x)=f*f^{*\nu}(x),\;\nu\ge1,$ while $f*g(x)=\int_0^x f(x-t)g(t)\,dt;$

3) Calculate the function $M(x)$ by the formula
\begin{equation}\label{MN}
M(x)=2N(x)-\int_0^x N^{*2}(t)\,dt, \quad 0<x<\pi.
\end{equation}}

The first part of Theorem~4 also was announced in \cite{But04}. Historically, it was the first result giving necessary and sufficient
conditions for solvability of an inverse spectral problems for an integro-differential operator. Equation (\ref{main_eq}) is called {\it
main nonlinear integral equation} of the inverse problem. In \cite{But07} (see also\cite{But06_MZ}), its global solvability was
established, i.e. for any function $w(x)\in L_2(0,\pi)$ equation (\ref{main_eq}) has a unique solution $N(x)\in L_{2,\pi},$ which has
played a crucial role in the proof of Theorem~4. Later on, developing the approach in \cite{But07} allowed researches to obtain global
solution also for other classes of integro-differential operators \cite{But10,BCh15, BS, BB17, BondBut-18, But18, But18-2, Bon18-1, Ign18,
Ign18-2, Bon18-2, Bon19-1, Bon19-2, Bon19-3,But19}. For different classes of operators, the corresponding main equations take different
forms, which makes it necessary to provide the proof of their solvability in each new case. In order to make it more convenient, in
\cite{ButMal18} a general approach has been developed for solving nonlinear equations of this type by introducing some abstract equation
and proving its global solvability. Moreover, in \cite{ButMal18} uniform stability of such nonlinear equations was established. In
\cite{ButTer19}, solvability of equation (\ref{main_eq}) was established in the class of entire functions of exponential type as soon as so
is the free term $w(x)$ and $w(0)=0.$ This fact along with the stability give an algorithm for solving equation (\ref{main_eq}) in
$L_{2,\pi},$ which can easily be implemented numerically.

The main result of the present paper is the following theorem, which gives uniform stability of Inverse Problem~1 in the weighted
$L_2$-norm as well as in the weighted uniform norm:
$$
\|f\|_{2,\pi}:=\|(\pi-x)f(x)\|_2, \quad\|f\|_{\infty,\pi}:=\|(\pi-x)f(x)\|_\infty,
$$
where we denoted $\|\cdot\|_2:=\|\cdot\|_{L_2(0,\pi)}$ and $\|\cdot\|_\infty:=\|\cdot\|_{L_\infty(0,\pi)},$ while the corresponding
distances between the spectra are determined as follows:
$$
\Lambda(\{\lambda_n\},\{\tilde\lambda_n\}):=\sqrt{\sum_{k=1}^\infty\frac{|\lambda_k-\tilde\lambda_k|^2}{k^2}}, \quad
\Lambda_1(\{\lambda_n\},\{\tilde\lambda_n\}):=\sum_{k=1}^\infty\frac{|\lambda_k-\tilde\lambda_k|}k,
$$
where $\{\lambda_k\}$ and $\{\tilde\lambda_k\}$ are the spectra of the problems ${\mathcal L}(M)$ and ${\mathcal L}(\tilde M).$

\medskip
{\bf Theorem 5. }{\it For any $r>0$ there exists $C_r>0$ such that
\begin{equation}\label{un_st}
\|M-\tilde M\|_{2,\pi}\le C_r\Lambda(\{\lambda_n\},\{\tilde\lambda_n\}), \quad \|M-\tilde M\|_{\infty,\pi}\le
C_r\Lambda_1(\{\lambda_n\},\{\tilde\lambda_n\})
\end{equation}
as soon as $\Lambda(\{\lambda_n\},\{n^2\})\le r$ and $\Lambda(\{\tilde\lambda_n\},\{n^2\})\le r.$}

\medskip
The second inequality in (\ref{un_st}) as well as other inequalities below with possibly infinite right-hand side mean that if it is
finite, then so is the left-hand side and the two are related as stated. Note that under the additional smoothness conditions on $M(x)$ and
$\tilde M(x)$ stated in the second part of Theorem~4, the value $\Lambda_1(\{\lambda_n\},\{\tilde\lambda_n\})$ is finite always when
$M(0)=\tilde M(0).$ The proof of Theorem~5 is based on Algorithm~1 and is a direct corollary of Lemmas~1--4 below, which give stability of
its steps 1)--3) in appropriate metrics.

The paper is organized as follows. In the next section, we provide Lemma~1, which gives stability of the first step in Algorithm~1. In
Section~3, we prove stability of steps~2) and~3) (Lemmas~2,~3 and~4, respectively). Throughout the paper, one and the same symbol $C_r$
denotes different positive constants in estimates, which depend only on $r.$
\\

{\large\bf 2. Stability of the characteristic function kernel}
\\

The function $\Delta(\lambda)$ involved in the first step of Algorithm~1 is called {\it characteristic function} of the problem ${\mathcal
L},$ whose eigenvalues coincide with zeros of $\Delta(\lambda)$ with account of multiplicity. Originally, it is determined as
$\Delta(\lambda)=S(\pi,\lambda),$ where $y=S(x,\lambda)$ is a solution of the equation in (\ref{eq1}) under the initial conditions
$S(0,\lambda)=0$ and $S'(0,\lambda)=1.$ In \cite{But07}, the representation
\begin{equation}\label{w}
\Delta(\lambda)=\frac{\sin\rho\pi}\rho +\int_0^\pi w(x)\frac{\sin\rho x}\rho\,dx, \quad w(x)\in L_2(0,\pi), \quad \rho^2=\lambda,
\end{equation}
was established, which, in turn, gives the asymptotics (\ref{eq3}) as well as the second formula in~(\ref{ser w}). Moreover, for any
sequence of complex numbers $\{\lambda_k\}$ of the form (\ref{eq3}), the function $\Delta(\lambda)$ constructed as the infinite product
in~(\ref{ser w}) has the form (\ref{w}) with a certain function $w(x)\in L_2(0,\pi),$ which, in turn, is determined by the Fourier series
in (\ref{ser w}) (see Lemma~3.3 in \cite{But07}).

Consider another sequence $\{\tilde\lambda_k\}$ of the form (\ref{eq3}) along with the corresponding functions
\begin{equation}\label{tilde_Delta}
\tilde \Delta(\lambda)=\pi\prod_{k=1}^\infty\frac{\tilde\lambda_k-\lambda}{k^2} =\frac{\sin\rho\pi}\rho +\int_0^\pi \tilde
w(x)\frac{\sin\rho x}\rho\,dx, \quad \tilde w(x)\in L_2(0,\pi).
\end{equation}
The following lemma gives uniform stabilities in $L_2$-metric and $L_\infty$-metric of recovering the kernel $w(x)$ from zeros
$\{\lambda_k\}$ of the function $\Delta(\lambda).$

\medskip
{\bf Lemma 1. }{\it For any $r>0,$ the following estimates hold:
\begin{equation}\label{est_w}
\|w-\tilde w\|_2\le C_r\Lambda(\{\lambda_n\},\{\tilde\lambda_n\}), \quad \|w-\tilde w\|_\infty \le
C_r\Lambda_1(\{\lambda_n\},\{\tilde\lambda_n\})
\end{equation}
as soon as $\Lambda(\{\lambda_n\},\{n^2\})\le r$ and $\Lambda(\{\tilde\lambda_n\},\{n^2\})\le r.$ Moreover,
$\Lambda_1(\{\lambda_n\},\{\tilde\lambda_n\})<\infty$ implies $w(x)-\tilde w(x) \in C[0,\pi].$}

\medskip
Before proceeding directly to the proof of Lemma~1, we provide several auxiliary assertions. First of all, we
prove the first estimate in (\ref{est_w}) in the particular case when $\tilde\lambda_k=k^2,\,k\ge1.$

\medskip
{\bf Proposition 1. }{\it For any $r>0,$ the estimate $\|w\|_2\le C_r\Lambda(\{\lambda_n\},\{n^2\})$ is
fulfilled as soon as $\Lambda(\{\lambda_n\},\{n^2\})\le r.$}

\medskip
{\it Proof.} According to (\ref{ser w}) and Parseval's equality, we calculate
$$
\|w\|_2=\sqrt{\frac2\pi\sum_{k=1}^\infty|k\Delta(k^2)|^2}, \quad k\Delta(k^2) =\pi k\prod_{j=1}^\infty
\frac{\lambda_j-k^2}{j^2}= a_k b_k\frac{\lambda_k-k^2}k, \quad a_k=\prod_{{j\ne k}\atop{j=1}}^\infty
\frac{\lambda_j-k^2}{j^2-k^2},
$$
$$
b_k=\pi\prod_{{j\ne k}\atop{j=1}}^\infty \frac{j^2-k^2}{j^2} =k^2\lim_{\rho\to
k}\frac{\sin\rho\pi}{\rho(k+\rho)(k-\rho)} =-\frac\pi2\lim_{\rho\to k}\cos\rho\pi =(-1)^{k+1}\frac\pi2.
$$
Thus, it remains to prove that $|a_k|\le C_r$ uniformly as $\Lambda(\{\lambda_n\},\{n^2\})\le r.$ For this purpose, we represent
$a_k=a_{k,1}a_{k,2},$ where
$$
a_{k,1}=\prod_{|j-k|\ge2r}\Big(1+\frac{\lambda_j-j^2}{j^2-k^2}\Big), \quad
a_{k,2}=\prod_{0<|j-k|<2r}\Big(1+\frac{\lambda_j-j^2}{j^2-k^2}\Big).
$$
Since
$$
\Big|\frac{\lambda_j-j^2}{j^2-k^2}\Big| =\frac{|\lambda_j-j^2|}j\frac{j}{|j-k|(j+k)}
\le\frac{\Lambda(\{\lambda_n\},\{n^2\})}{2r} \le\frac12, \quad |j-k|\ge2r,
$$
we get the estimate
$$
|a_{k,1}|=\Big|\exp\Big(\sum_{|j-k|\ge2r}\ln\Big(1+\frac{\lambda_j-j^2}{j^2-k^2}\Big)\Big)\Big|\le
\exp\Big(2\sum_{|j-k|\ge2r}\Big|\frac{\lambda_j-j^2}{j^2-k^2}\Big|\Big),
$$
where the Cauchy--Bunyakovsky--Schwarz inequality implies
$$
\sum_{|j-k|\ge2r}\Big|\frac{\lambda_j-j^2}{j^2-k^2}\Big| \le \sqrt{\sum_{|j-k|\ge2r}
\frac{|\lambda_j-j^2|^2}{j^2} \sum_{|j-k|\ge2r}\frac{j^2}{(j^2-k^2)^2}}\le r\sqrt{ \sum_{{j\ne
k}\atop{j=1}}^\infty\frac{j^2}{(j^2-k^2)^2}},
$$
while
$$
\sum_{{j\ne k}\atop{j=1}}^\infty\frac{j^2}{(j^2-k^2)^2} =\sum_{{j\ne
k}\atop{j=1}}^\infty\frac{j^2}{(j-k)^2(j+k)^2} <\sum_{j=1}^{k-1}\frac1{(k-j)^2}
+\sum_{j=k+1}^\infty\frac1{(j-k)^2} <2\sum_{j=1}^\infty\frac1{j^2} =\frac{\pi^2}3.
$$
Finally, we get
$$
|a_{k,2}|\le\prod_{0<|j-k|<2r}\Big(1+\frac{rj}{|j-k|(j+k)}\Big)
 <(1+r)^{4r-2},
$$
which finishes the proof. $\hfill\Box$

\medskip
In the general case, the proof begins with the following assertion.

\medskip
{\bf Proposition 2. }{\it There exists a choice of $\{\varkappa_k\}$ in (\ref{eq3}) such that for any $r>0$
the estimate $|\varkappa_k|\le r$ holds for all $k\in{\mathbb N}$ as soon as
$\Lambda(\{\lambda_n\},\{n^2\})\le r.$}

\medskip
{\it Proof.} Putting $\varepsilon_k:=(\lambda_k-k^2)/k,$ we arrive at $|\varepsilon_k|\le \Lambda(\{\lambda_n\},\{n^2\})\le r.$ Further, we
have $\lambda_k=k^2+k\varepsilon_k =(k+\varkappa_k)^2,$ where $\varkappa_k =k(\sqrt{1+\varepsilon_k/k}-1)$ and ${\rm
Re}\sqrt{\,\cdot\,}\ge0.$ Hence, $|\varkappa_k|=|\varepsilon_k|/|\sqrt{1+\varepsilon_k/k} +1| \le|\varepsilon_k|\le r.$ $\hfill\Box$

\medskip
In what follows, without loss of generality we assume that $r\in{\mathbb N}.$

\medskip
{\bf Proposition 3. }{\it For any $r\in{\mathbb N},$ the estimate $|\lambda_j-k^2|\ge 4jr$ holds as soon as $|j-k|\ge 6r$ and
$\Lambda(\{\lambda_n\},\{n^2\})\le r.$}

\medskip
{\it Proof.} According to (\ref{eq3}) and Proposition~2 we have
$|\lambda_j-k^2|=|j-k+\varkappa_j||j+k+\varkappa_j| \ge(|j-k|-r)(j+k-r).$ Thus, it is sufficient to prove
that
\begin{equation}\label{ineq-1}
(|j-k|-r)(j+k-r)\ge4jr.
\end{equation}
Let $j\ge k.$ Then (\ref{ineq-1}) is equivalent to the inequality $j^2-6jr-k^2+r^2\ge0,$ which is, obviously,
fulfilled for $j-k\ge 6r,$ by virtue of the inequalities $j\ge k+6r > 3r +\sqrt{8r^2 +k^2}.$

Further, let $k>j.$ Then (\ref{ineq-1}) is equivalent to the inequality $j^2+4jr-(k-r)^2\le0,$ which, in
turn, holds for $k-j\ge 5r,$ because
$$
j\le k-5r =-2r+\sqrt{4r^2-4r(k-r) +(k-r)^2} <-2r+\sqrt{4r^2 +(k-r)^2}.
$$
Thus, for $|j-k|\ge 6r$ inequality (\ref{ineq-1}) is proven. $\hfill\Box$

\medskip
For $r,\,k\in{\mathbb N}$ we introduce the sets
$$
\Omega_r(k):=\{j:|j-k|<6r,j\in{\mathbb N}\}, \;\; \Omega_r'(k):=\Omega_r(k)\setminus\{k\}, \;\; \Theta_r(k):=\{j:|j-k|\ge 6r,j\in{\mathbb
N}\}.
$$
Obviously, $\Omega_r(k)\cup\Theta_r(k)\equiv{\mathbb N}.$ We also put $\alpha_{r,k}:=\#\Omega_r(k)=\min\{k,6r\}+6r-1\le 12r-1.$

\medskip
Denote
$$
\sigma_{r,k}(\lambda):=\prod_{j\in\Omega_r(k)}\frac{\lambda_j-\lambda}{j^2}, \quad
\tilde\sigma_{r,k}(\lambda):=\prod_{j\in\Omega_r(k)}\frac{\tilde\lambda_j-\lambda}{j^2}.
$$

\medskip
{\bf Proposition 4. }{\it For any $r\in{\mathbb N},$ the estimates
\begin{equation}\label{theta_est}
|\sigma_{r,k}(k^2)|\le \frac{C_r}{k^{\alpha_{r,k}}}\frac{|\lambda_k-k^2|}k, \quad |\sigma_{r,k}(k^2)
-\tilde\sigma_{r,k}(k^2)|\le\frac{C_r}{k^{\alpha_{r,k}}}\sum_{j\in\Omega_r(k)}\frac{|\lambda_j-\tilde\lambda_j|}j, \quad k\in{\mathbb N},
\end{equation}
are fulfilled as soon as $\Lambda(\{\lambda_n\},\{n^2\})\le r$ and $\Lambda(\{\tilde\lambda_n\},\{n^2\})\le r.$}

\medskip
{\it Proof.} We have
\begin{equation}\label{theta_est1}
\sigma_{r,k}(k^2)=\frac{\lambda_k-k^2}{k^2}\prod_{j\in\Omega_r'(k)}\frac{\lambda_j-k^2}{j^2}.
\end{equation}
Since $j\in\Omega_r(k)$ is equivalent to the inequalities $\max\{0,k-6r\}< j< k+6r,$ we have
\begin{equation}\label{theta_est2}
|\lambda_j-k^2| =|j^2-k^2+2j\varkappa_j +\varkappa_j^2| < 6r(2k+6r) +2(k+6r)r +r^2 \le C_rk, \quad
j\in\Omega_r(k), \;\; k\in{\mathbb N},
\end{equation}
and
\begin{equation}\label{theta_est3}
\frac1j\le\frac1{\max\{0,k-6r\}+1} \le\frac{6r}{k}, \quad j\in\Omega_r(k), \;\; k\in{\mathbb N}.
\end{equation}
Substituting estimates (\ref{theta_est2}) and (\ref{theta_est3}) into (\ref{theta_est1}), we get
$$
|\sigma_{r,k}(k^2)|\le \frac{|\lambda_k-k^2|}{k^2} \Big(\frac{36r^2C_r}k\Big)^{\alpha_{r,k}-1},
$$
which coincides with the first estimate in (\ref{theta_est}). Further, it is easy to check that
\begin{equation}\label{theta_est4}
\sigma_{r,k}(k^2) -\tilde\sigma_{r,k}(k^2)
=\sum_{j\in\Omega_r(k)}\tilde\sigma_{r,k,j}(k^2)\frac{\lambda_j-\tilde\lambda_j}{j^2}\sigma_{r,k,j}(k^2),
\end{equation}
where
$$
\sigma_{r,k,j}(\lambda)=\prod_{\nu=j+1}^{k+6r-1}\frac{\lambda_\nu-\lambda}{\nu^2}, \quad
\tilde\sigma_{r,k,j}(\lambda)=\prod_{\nu=\max\{0,k-6r\}+1}^{j-1}\frac{\tilde\lambda_\nu-\lambda}{\nu^2}.
$$
According to (\ref{theta_est2}) and (\ref{theta_est3}), we have
$$
|\sigma_{r,k,j}(k^2)\tilde\sigma_{r,k,j}(k^2)| \le\Big(\frac{36r^2C_r}k\Big)^{\alpha_{r,k}-1},
$$
which along with (\ref{theta_est3}) and (\ref{theta_est4}) give the second estimate in (\ref{theta_est}). $\hfill\Box$

\medskip
Denote
\begin{equation}\label{Delta_k}
\Delta_k(\lambda):=\frac{\Delta(\lambda)}{\sigma_{r,k}(\lambda)} =\pi\prod_{j\in\Theta_r(k)}\frac{\lambda_j-\lambda}{j^2}, \quad
\tilde\Delta_k(\lambda):=\frac{\tilde\Delta(\lambda)}{\tilde\sigma_{r,k}(\lambda)}
=\pi\prod_{j\in\Theta_r(k)}\frac{\tilde\lambda_j-\lambda}{j^2}.
\end{equation}

\medskip
{\bf Proposition 5. }{\it For any $r\in{\mathbb N},$ the estimate
\begin{equation}\label{Delta_k-0}
|\Delta_k(k^2)|\le C_r k^{\alpha_{r,k}-1}, \quad k\in{\mathbb N},
\end{equation}
holds as soon as $\Lambda(\{\lambda_n\},\{n^2\})\le r.$}

\medskip
{\it Proof.} Put $\varepsilon_r:= 8r+1.$ Then, in particular, according to (\ref{eq3}) and Proposition~2, we have
$|\rho_k-k|<\varepsilon_r$ for all $k\in{\mathbb N}.$ Hence, the maximum modulus principle gives
\begin{equation}\label{Delta_k-1}
|\Delta_k(k^2)| <\max_{|\rho-\rho_k|=\varepsilon_r}\Big|\frac{\Delta(\rho^2)}{\sigma_{r,k}(\rho^2)}\Big|
=\max_{|\rho-\rho_k|=\varepsilon_r}\Big|\Delta(\rho^2)\prod_{j\in\Omega_r(k)}\frac{j^2}{\lambda_j-\rho^2}\Big|.
\end{equation}
By virtue of representation (\ref{w}) and Proposition 1, we have the estimate
$$
|\Delta(\rho^2)|\le\frac{A_{r,C}}{|\rho|+r+\varepsilon_r}, \quad |{\rm Im}\rho|\le C,
$$
where $A_{r,C}$ depends only on $r$ and $C.$ The latter estimate holds also if $|\rho-\rho_k|=\varepsilon_r$ for any $k\in{\mathbb N},$
because in this case we have $|{\rm Im}\rho|\le|{\rm Im}\rho_k| +|{\rm Im}(\rho-\rho_k)|\le C:=r+\varepsilon_r.$ Furthermore, according to
Proposition~2, for $|\rho-\rho_k|=\varepsilon_r$ we have the estimate
$$
\frac1{|\rho|+r+\varepsilon_r} \le\frac1{|\rho_k|-|\rho-\rho_k|+r+\varepsilon_r} =\frac1{|\rho_k|+r}\le\frac1k, \quad k\in{\mathbb N}.
$$
Thus, we get the estimate
\begin{equation}\label{Delta_k-2}
|\Delta(\rho^2)|\le \frac{C_r}k, \quad |\rho-\rho_k|=\varepsilon_r, \quad k\in{\mathbb N}.
\end{equation}
Further, we have
\begin{equation}\label{Delta_k-3}
\prod_{j\in\Omega_r(k)}j^2 <(k+6r)^{2\alpha_{r,k}}\le C_r k^{2\alpha_{r,k}}, \quad k\in{\mathbb N}.
\end{equation}
Moreover, if $|\rho-\rho_k|=\varepsilon_r$ and $j\in\Omega_r(k),$ then we also have
$$
|\rho-\rho_j|\ge\varepsilon_r -|\rho_j-\rho_k| \ge\varepsilon_r -|j-k|-|\varkappa_j|-|\varkappa_k| >\varepsilon_r-8r=1, \quad k\ge1,
$$
$$
|\rho+\rho_j|\ge|\rho_j+\rho_k|-\varepsilon_r \ge j+k-2r -\varepsilon_r >2k-8r -\varepsilon_r\qquad\qquad\qquad\quad
$$
$$
\qquad\qquad\qquad=2k-2\varepsilon_r+1=\Big(2-\frac{2\varepsilon_r-1}k\Big)k\ge \Big(2-\frac{2\varepsilon_r-1}{\varepsilon_r}\Big)k
=\frac{k}{\varepsilon_r}, \quad k\ge\varepsilon_r.
$$
Hence, we have $|\lambda_j-\rho^2|\ge k/\varepsilon_r$ as soon as $|\rho-\rho_k|=\varepsilon_r,\;j\in\Omega_r(k)$ and $k\ge\varepsilon_r,$
which along with (\ref{Delta_k-1})--(\ref{Delta_k-3}) give (\ref{Delta_k-0}) for $k\ge\varepsilon_r.$ Further, for
$k=\overline{1,\varepsilon_r-1}$ and $j\in\Omega_r(k)$ we have
$$
j<k+6r\le14r, \quad |\lambda_j-\rho^2| =|\rho-\rho_j||\rho+\rho_j| \ge(16r-j-r)^2=r^2\quad {\rm for} \quad |\rho|=16r.
$$
Hence, for $k<\varepsilon_r,$ estimate (\ref{Delta_k-0}) follows from the following rough estimate:
$$
|\Delta_k(k^2)| <\max_{|\rho|=16r}\Big|\Delta(\rho^2)\prod_{j\in\Omega_r(k)}\frac{j^2}{\lambda_j-\rho^2}\Big|
\le14^{24r-2}\max_{|\rho|=16r}|\Delta(\rho^2)|\le C_r, \quad k=\overline{1,\varepsilon_r-1}.
$$
Thus, we arrive at (\ref{Delta_k-0}) for all $k\in{\mathbb N}.$  $\hfill\Box$

\medskip
Denote
$$
\theta_k:=\sum_{j\in\Theta_r(k)}\Big|\frac{\lambda_j-\tilde\lambda_j}{\lambda_j-k^2}\Big|, \quad k\in{\mathbb N}.
$$
The following proposition gives an estimate for the sequence $\{\theta_k\}$ in the $l_\infty$-norm.

\medskip
{\bf Proposition 6. }{\it For any $r\in{\mathbb N},$ the estimate
$$
\sup_{k\in{\mathbb N}}\theta_k \le C_r\Lambda(\{\lambda_n\},\{\tilde\lambda_n\})
$$
is fulfilled as soon as $\Lambda(\{\lambda_n\},\{n^2\})\le r.$}

\medskip
{\it Proof.} The Cauchy--Bunyakovsky--Schwarz inequality gives
$$
\theta_k \le\sqrt{\sum_{j\in\Theta_r(k)}\frac{j^2}{|\lambda_j-k^2|^2} \sum_{j\in\Theta_r(k)}\frac{|\lambda_j-\tilde\lambda_j|^2}{j^2}} \le
\alpha_k\Lambda(\{\lambda_n\},\{\tilde\lambda_n\}), \quad \alpha_k=\sqrt{\sum_{j\in\Theta_r(k)}\frac{j^2}{|\lambda_j-k^2|^2}}.
$$
Thus, it remains to prove that $\alpha_k\le C_r$ for $k\in{\mathbb N}.$ We have $\alpha_k^2=\eta_k+\zeta_k,$ where $\eta_k=0$ for $k\le6r$
and
$$
\eta_k=\sum_{j=1}^{k-6r}\frac{j^2}{|\lambda_j-k^2|^2} =\sum_{j=1}^{k-6r}\frac{j^2}{(\rho_j-k)^2(\rho_j+k)^2}
\le\sum_{j=1}^{k-6r}\frac{(k-6r)^2}{(k-j-r)^2(k+j-r)^2} \qquad\qquad
$$
$$
\qquad\qquad\qquad\qquad\qquad\qquad\qquad\qquad <\sum_{j=1}^{k-6r}\frac1{(k-j-r)^2} =\sum_{j=5r}^{k-r-1}\frac1{j^2} <\frac{\pi^2}6, \quad
k>6r,
$$
while
$$
\zeta_k =\sum_{j=k+6r}^\infty\frac{j^2}{|\lambda_j-k^2|^2} \le\sum_{j=k+6r}^\infty\frac{j^2}{(j-k-r)^2(j+k-r)^2}, \quad k\in{\mathbb N}.
$$
Since
$$
\frac{j}{j+k-r}\le \frac1{1-\frac{r-1}j}\le \frac1{1-\frac{r-1}{1+6r}} =\frac{1+6r}{2+5r}<\frac65, \quad j\ge 1+6r,
$$
we arrive at
$$
\zeta_k \le\frac{36}{25}\sum_{j=k+6r}^\infty\frac1{(j-k-r)^2} =\frac{36}{25}\sum_{j=1}^\infty\frac1{(j+5r-1)^2}<\frac{6\pi^2}{25},
$$
which finishes the proof.  $\hfill\Box$

\medskip
Finally, we estimate the sequence $\{\theta_k\}$ also in the $l_2$-norm.

\medskip
{\bf Proposition 7. }{\it For any $r>0,$ the estimate
$$
\sqrt{\sum_{k=1}^\infty\theta_k^2} \le C_r\Lambda_1(\{\lambda_n\},\{\tilde\lambda_n\})
$$
is fulfilled as soon as $\Lambda(\{\lambda_n\},\{n^2\})\le r.$}

\medskip
{\it Proof.} For $k\in{\mathbb N},$ putting
$$
\beta_{k,j}:=\Big|\frac{\lambda_j-\tilde\lambda_j}{\lambda_j-k^2}\Big|, \quad j\in\Theta_r(k), \qquad \beta_{k,j}:=0, \quad
j\in\Omega_r(k),
$$
and using the generalized Minkowski inequality, we get
$$
\sqrt{\sum_{k=1}^\infty\theta_k^2}=\sqrt{\sum_{k=1}^\infty\Big(\sum_{j\in\Theta_r(k)}\Big|
\frac{\lambda_j-\tilde\lambda_j}{\lambda_j-k^2}\Big|\Big)^2}= \sqrt{\sum_{k=1}^\infty\Big(\sum_{j=1}^\infty\beta_{k,j}\Big)^2} \le
\sum_{j=1}^\infty\sqrt{\sum_{k=1}^\infty\beta_{k,j}^2} \qquad\qquad
$$
$$
\qquad\qquad=\sum_{j=1}^\infty\sqrt{\sum_{k\in\Theta_r(j)}\Big|\frac{\lambda_j-\tilde\lambda_j}{\lambda_j-k^2}\Big|^2} =\sum_{j=1}^\infty
\gamma_j\frac{|\lambda_j-\tilde\lambda_j|}j, \quad \gamma_j=\sqrt{\sum_{k\in\Theta_r(j)}\frac{j^2}{|\lambda_j-k^2|^2}}.
$$
Thus, it is sufficient to prove that $\gamma_j\le C_r$ for $j\in{\mathbb N}.$ We have $\gamma_j^2=\mu_j+\xi_j,$ where
$$
\mu_j =\sum_{k=j+6r}^\infty\frac{j^2}{|\lambda_j-k^2|^2} \le \sum_{k=j+6r}^\infty\frac{j^2}{(|\rho_j|+k)^2(|\rho_j|-k)^2} \le
\frac{j^2}{(|\rho_j|+j+6r)^2}\sum_{k=j+6r}^\infty\frac1{(|\rho_j|-k)^2}
$$
$$
<\sum_{k=j+6r}^\infty\frac1{(k-j-r)^2} =\sum_{k=1}^\infty\frac1{(k+5r-1)^2} <\frac{\pi^2}6,
$$
while $\xi_j=0$ for $j\ge6r$ and
$$
\xi_j=\sum_{k=1}^{j-6r}\frac{j^2}{|\lambda_j-k^2|^2} \le \sum_{k=1}^{j-6r}\frac{j^2}{(|\lambda_j|-k^2)^2}
=\sum_{k=1}^{j-6r}\frac{j^2}{(|\rho_j|+k)^2(|\rho_j|-k)^2} \qquad\qquad\qquad\;
$$
$$
\qquad \le\frac{j^2}{|\lambda_j|+1}\sum_{k=1}^{j-6r}\frac1{(j-k-r)^2} =\frac{j^2}{|\lambda_j|+1}\sum_{k=5r}^{j-r-1}\frac1{k^2}
<\frac{\pi^2}6\frac{j^2}{|\lambda_j|+1}, \quad j>6r.
$$
It remains to note that $(|\lambda_j|+1)^{-1}j^2\le r^2$ for $j\le r$ and
$$
\frac{j^2}{|\lambda_j|+1}< \frac{j^2}{(j-r)^2} =  \frac1{(1-\frac{r}j)^2} \le \frac1{(1-\frac{r}{r+1})^2} =(r+1)^2
$$
for $j\ge r+1.$ $\hfill\Box$

\medskip
Now we are in position to give the proof of Lemma~1.

\medskip
{\it Proof of Lemma 1}. By virtue of (\ref{ser w}) and (\ref{tilde_Delta}) along with Parseval's equality, we have
\begin{equation}\label{1+}
\hat w(x)=\frac2\pi\sum_{k=1}^\infty k\hat\Delta(k^2)\sin kx , \quad \|\hat w\|_2= \sqrt{\frac2\pi\sum_{k=1}^\infty |k\hat\Delta(k^2)|^2},
\end{equation}
where $\hat w(x) =w(x)-\tilde w(x)$ and $\hat\Delta(\lambda) =\Delta(\lambda)-\tilde \Delta(\lambda).$ According to (\ref{Delta_k}), we
arrive at
\begin{equation}\label{hat_Delta}
\hat\Delta(k^2) =\Delta_k(k^2)\Big(\sigma_{r,k}(k^2) -\tilde\sigma_{r,k}(k^2)
+\Big(1-\frac{\tilde\Delta_k(k^2)}{\Delta_k(k^2)}\Big)\tilde\sigma_{r,k}(k^2)\Big),
\end{equation}
where
$$
\Big|1-\frac{\tilde\Delta_k(k^2)}{\Delta_k(k^2)}\Big| =\Big|1-\prod_{j\in\Theta_r(k)}\frac{\tilde\lambda_j-k^2}{\lambda_j-k^2}\Big|
=\Big|1-\exp\Big(\sum_{j\in\Theta_r(k)}\ln\Big(1-\frac{\lambda_j-\tilde\lambda_j}{\lambda_j-k^2}\Big)\Big)\Big|,
$$
where, in turn, by virtue of Proposition~3, we have
$$
\Big|\frac{\lambda_j-\tilde\lambda_j}{\lambda_j-k^2}\Big| =\frac{|\lambda_j-\tilde\lambda_j|}j\frac{j}{|\lambda_j-k^2|} \le
\frac{\Lambda(\{\lambda_n\},\{n^2\}) +\Lambda(\{\tilde\lambda_n\},\{n^2\})}{4r}\le\frac12, \quad j\in\Theta_r(k).
$$
Thus, we get
$$
\Big|1-\frac{\tilde\Delta_k(k^2)}{\Delta_k(k^2)}\Big|
\le\sum_{\nu=1}^\infty\frac{2^\nu}{\nu!}\Big(\sum_{j\in\Theta_r(k)}\Big|\frac{\lambda_j-\tilde\lambda_j}{\lambda_j-k^2}\Big|\Big)^\nu
=\sum_{\nu=1}^\infty\frac{(2\theta_k)^\nu}{\nu!} \le2\theta_k\exp(2\theta_k),
$$
which along with (\ref{hat_Delta}) and Propositions~4 and~5 give
\begin{equation}\label{hat_Delta-1}
k|\hat\Delta(k^2)| \le C_r\sum_{j\in\Omega_r(k)}\frac{|\lambda_j-\tilde\lambda_j|}j +C_r\theta_k\exp(2\theta_k)\frac{|\lambda_k-k^2|}k.
\end{equation}
Since
\begin{equation}\label{inequalities}
\Big(\sum_{k=1}^n a_k\Big)^2\le n\sum_{k=1}^n a_k^2, \qquad \sum_{k=1}^\infty\sum_{j\in\Omega_r(k)} a_{k,j} =\sum_{k,j\in{\mathbb N},\,
|k-j|<6r} a_{k,j} =\sum_{j=1}^\infty\sum_{k\in\Omega_r(j)} a_{k,j}
\end{equation}
for, in particular, any non-negative summands, we have
$$
\sqrt{\sum_{k=1}^\infty\Big(\sum_{j\in\Omega_r(k)}\frac{|\lambda_j-\tilde\lambda_j|}j\Big)^2} \le
\sqrt{\sum_{k=1}^\infty\alpha_{r,k}\sum_{j\in\Omega_r(k)}\frac{|\lambda_j-\tilde\lambda_j|^2}{j^2}} \qquad\qquad\qquad\qquad\qquad\qquad
$$
$$
\quad\qquad\qquad\qquad\qquad\qquad =\sqrt{\sum_{j=1}^\infty\frac{|\lambda_j-\tilde\lambda_j|^2}{j^2}\sum_{k\in\Omega_r(j)}\alpha_{r,k}}
\le(12r-1)\Lambda(\{\lambda_n\},\{\tilde\lambda_n\}),
$$
which along with (\ref{hat_Delta-1}) and Proposition 6 give
$$
\sqrt{\sum_{k=1}^\infty |k\hat\Delta(k^2)|^2} \le C_r\Big(12r-1 +C_r\Lambda(\{\lambda_n\},\{n^2\})\Big)
\Lambda(\{\lambda_n\},\{\tilde\lambda_n\}).
$$
Taking the second equality in (\ref{1+}) into account, we arrive at the first estimate in (\ref{est_w}).

Further, by virtue of the first equality in (\ref{1+}) along with (\ref{hat_Delta-1}) and the Cauchy--Bunyakovsky--Schwarz inequality, we
get
\begin{equation}\label{2+}
\|\hat w\|_\infty \le\frac2\pi\sum_{k=1}^\infty k|\hat\Delta(k^2)| \le
C_r\sum_{k=1}^\infty\sum_{j\in\Omega_r(k)}\frac{|\lambda_j-\tilde\lambda_j|}j +C_r
\Lambda(\{\lambda_n\},\{n^2\})\sqrt{\sum_{k=1}^\infty\theta_k^2}.
\end{equation}
According to (\ref{inequalities}) and Proposition~7, we obtain
$$
\sum_{k=1}^\infty k|\hat\Delta(k^2)| \le C_r\sum_{j=1}^\infty\sum_{k\in\Omega_r(j)}\frac{|\lambda_j-\tilde\lambda_j|}j +C_r
\Lambda_1(\{\lambda_n\},\{\tilde\lambda_n\}) \le 12rC_r\Lambda_1(\{\lambda_n\},\{\tilde\lambda_n\}),
$$
which along with (\ref{2+}) give the second estimate in (\ref{est_w}). Finally, note that
$\Lambda_1(\{\lambda_n\},\{\tilde\lambda_n\})<\infty$ implies the uniform convergence of the series in (\ref{1+}), which gives $\hat
w(x)\in C[0,\pi].$  $\hfill\Box$
\\

{\large\bf 3. Stability of steps 2) and 3) in Algorithm 1}
\\

Denote $h(x):=(\pi-x)N(x)$ and consider the nonlinear operator
$$
{\mathcal D}h(x):=\sum_{\nu=2}^\infty \frac{(\pi-x)^\nu}{\nu!} N^{*\nu}(x).
$$
Then the main equation (\ref{main_eq}) takes the form $w(\pi-x)=h(x)+{\mathcal D}h(x),$ while the operator ${\mathcal D}$ belongs to the
class ${\mathcal E}'_{\pi,1}$ in \cite{ButMal18} and satisfies Condition ${\mathcal A}$ therein. Thus, the following lemma, giving uniform
stability of equation (\ref{main_eq}) in the metric of $L_{2,\pi},$ follows from Corollary~1 in \cite{ButMal18}.

\medskip
{\bf Lemma 2. }{\it For any $r>0$ the estimate $\|N-\tilde N\|_{2,\pi}\le C_r\|w-\tilde w\|_2$ holds as soon as $\|w\|_2\le r$ and
$\|\tilde w\|_2\le r,$ where $N(x)$ is the solution of equation (\ref{main_eq}), while $\tilde N(x)$ is the one of the equation
\begin{equation}\label{main_eq_ti}
\tilde w(\pi-x)=\sum_{\nu=1}^\infty \frac{(\pi-x)^\nu}{\nu!} \tilde N^{*\nu}(x).
\end{equation}}

\medskip
The following lemma gives uniform stability of equation (\ref{main_eq}) in the weighted uniform norm.

\medskip
{\bf Lemma 3. }{\it For any $r>0$ the estimate $\|N-\tilde N\|_{\infty,\pi}\le C_r\|w-\tilde w\|_\infty$ holds as soon as $\|w\|_2\le r$
and $\|\tilde w\|_2\le r,$ where $N(x)$ and $\tilde N(x)$ are solutions of equations (\ref{main_eq}) and~(\ref{main_eq_ti}). }

\medskip {\it Proof.} Subtracting (\ref{main_eq_ti}) from (\ref{main_eq}), we get
$$
\hat w(\pi-x)=\hat h(x)+\int_0^x K(x,t)\hat h(t)\,dt,
$$
where $\hat w(x)=w(x)-\tilde w(x),$ $\hat h(x)=h(x)-\tilde h(x),$ $\tilde h(x)=(\pi-x)\tilde N(x)$ and
$$
K(x,t)= \frac1{\pi-t}\sum_{\nu=2}^\infty \frac{(\pi-x)^\nu}{\nu!} f_\nu(x-t), \quad f_\nu=N^{*(\nu-1)} +\sum_{j=1}^{\nu-2}N^{*j}*\tilde
N^{*(\nu-1-j)} +\tilde N^{*(\nu-1)}.
$$
By virtue of Lemma~2, there exists $r_0>0$ depending only on $r$ such that $\|w\|_2\le r$ implies $\|h(x)\|_2=\|N\|_{2,\pi}\le r_0.$ Let
$R(x,t)$ be the resolvent kernel for $K(x,t),$ i.e.
\begin{equation}\label{5-2}
\hat h(x)=\hat w(\pi-x)+\int_0^x R(x,t)\hat w(\pi-t)\,dt, \quad R(x,t)=-K(x,t) -\int_t^x K(x,\tau)R(\tau,t)\,d\tau.
\end{equation}
Thus, since $\|N-\tilde N\|_{\infty,\pi}=\|\hat h\|_\infty,$ it is sufficient to prove the estimate
\begin{equation}\label{5-3}
A_R \le C_r, \quad {\rm where} \quad A_G:={\rm ess}\sup_{0<x<\pi}\sqrt{\int_0^x |G(x,t)|^2\,dt}.
\end{equation}
Indeed, the second equality in (\ref{5-2}) implies the estimate
$$
|R(x,t)|\le|K(x,t)|+ A_K \sqrt{\int_t^x |R(\tau,t)|^2\,d\tau},
$$
which implies $A_R\le (1+R)A_K,$ where $R=\|R(x,t)\|_{L_2((0,\pi)^2)}.$ Moreover, by virtue of Lemma~1 in \cite{ButMal18}, we have $R\le
F(K),$ where $F(x)=x+ \sum_{n=0}^\infty \frac{x^{n+2}}{\sqrt{n!}}$ and $K=\|K(x,t)\|_{L_2((0,\pi)^2)}.$ Thus, it remains to note that $K\le
\sqrt\pi A_K$ and
$$
(\pi-x)^{\nu-1}|f_\nu(x-t)|\le \nu(\pi-x)^{\nu-1}N_1^{*(\nu-1)}(x-t)\le \nu h_1^{*(\nu-1)}(x-t),
$$
where $N_1(x)=|N(x)|+|\tilde N(x)|$ and $h_1(x)=(\pi-x)N_1(x)=|h(x)|+|\tilde h(x)|.$ Hence, we have $\|h_1\|_2\le2r_0$ and arrive at the
estimates
$$
|K(x,t)|\le \sum_{\nu=2}^\infty \frac{h_1^{*(\nu-1)}(x-t)}{(\nu-1)!}, \quad A_K\le \sum_{\nu=1}^\infty \frac{\|h_1^{*\nu}\|_2}{\nu!}
\le\sum_{\nu=1}^\infty \frac{\pi^\frac{\nu-1}2\|h_1\|_2^\nu}{\nu!} \le\frac{\exp(2r_0\sqrt\pi)-1}{\sqrt\pi},
$$
which finish the proof of (\ref{5-3}). $\hfill\Box$

\medskip
The following lemma gives the uniform stabilities of the third step in Algorithm~1.

\medskip
{\bf Lemma 4. }{\it For any $r>0$ the following estimates hold:
\begin{equation}\label{hat_M}
\|M-\tilde M\|_{2,\pi}\le C_r\|N-\tilde N\|_{2,\pi}, \quad \|M-\tilde M\|_{\infty,\pi}\le C_r\|N-\tilde N\|_{\infty,\pi}
\end{equation}
as soon as $\|N\|_{2,\pi}\le r$ and $\|\tilde N\|_{2,\pi}\le r,$ where the function $M(x)$ is determined by $N(x)$ via formula~(\ref{MN}),
while $\tilde M(x)$ is determined by $\tilde N(x)$ via the analogous formula
\begin{equation}\label{tildeMN}
\tilde M(x)=2\tilde N(x)-\int_0^x \tilde N^{*2}(t)\,dt.
\end{equation}}

{\it Proof.} For briefness, we denote $\hat N(x):=N(x)-\tilde N(x)$ and $\hat M(x):=M(x)-\tilde M(x).$ Then, subtracting (\ref{tildeMN})
from (\ref{MN}), we get
\begin{equation}\label{hat_M-1}
\hat M(x)=2\hat N(x)-1*( N^{*2}-\tilde N^{*2})(x)= 2\hat N(x)-1*f_2*\hat N(x),
\end{equation}
where $f_2(x)=N(x)+\tilde N(x)$ and, hence, $\|f_2\|_{2,\pi}\le 2r.$ Thus, we arrive at the estimate
$$
(\pi-x)|1*f_2*\hat N(x)|\le \int_0^x|f_2(t)|\,dt \int_0^{x-t} (\pi-\tau)|\hat N(\tau)|\,d\tau \qquad\qquad\qquad \qquad\qquad\qquad\qquad
$$
\begin{equation}\label{hat_M-2}
\qquad\qquad\le \int_0^x |f_2(t)|\sqrt{\int_0^{x-t}(\pi-\tau)^2|\hat N(\tau)|^2\,d\tau}\sqrt{x-t}\,dt \le \|\hat N\|_{2,\pi}
\int_0^x|f_2(t)|\sqrt{x-t}\,dt.
\end{equation}
We also get
$$
\int_0^\pi\Big(\int_0^x|f_2(t)|\sqrt{x-t}\,dt\Big)^2dx \le \int_0^\pi x\,dx\int_0^x|f_2(t)|^2(x-t)\,dt \qquad\qquad\qquad\qquad\qquad\qquad
$$
$$
\qquad\quad=\int_0^\pi |f_2(x)|^2\,dx\int_x^\pi t(t-x)\,dt\le\pi\int_0^\pi (\pi-x)|f_2(x)|^2\,dx\int_x^\pi dt =\pi\|f_2\|_{2,\pi}^2\le 4\pi
r^2,
$$
which along with (\ref{hat_M-1}) and (\ref{hat_M-2}) give the first estimate in (\ref{hat_M}) with $C_r=2+2r\sqrt\pi.$

Further, by virtue of (\ref{hat_M-1}), we get the estimate
$$
(\pi-x)|\hat M(x)|\le2(\pi-x)|\hat N(x)|+\|\hat N\|_{\infty,\pi}\int_0^x (x-t)|f_2(t)|\,dt,
$$
which gives $\|\hat M\|_{\infty,\pi} \le 2(1+r\sqrt\pi)\|\hat N\|_{\infty,\pi}$ and finishes the proof. $\hfill\Box$
\\

{\bf Funding.} This work was supported by Grant 20-31-70005 of the Russian Foundation for Basic Research.

\end{document}